\documentclass[onecolumn,final,a4paper,sort&compress]{elsarticle}

\usepackage[utf8]{inputenc}
\usepackage{amsmath,amssymb,amsfonts,latexsym,stmaryrd}
\usepackage{algorithmic}
\usepackage{booktabs}

\usepackage[colorlinks=true,breaklinks=true,linkcolor=lightblue,citecolor=lightblue]{hyperref}
\usepackage{stackengine}
\usepackage{mathrsfs}
\usepackage{rotating} 
\usepackage{showkeys}
\usepackage{amssymb,amsmath,bbm,amsfonts,latexsym,subfigure}
\usepackage{amsthm}
\usepackage{graphicx,float,epsfig,color,fancyhdr}
\usepackage{amsopn}

\usepackage[utf8]{inputenc}
\usepackage{graphicx}
\usepackage{wrapfig}

\usepackage{amsfonts}
\usepackage{amsmath}
\numberwithin{equation}{section}
\usepackage{mathtools}

\usepackage[ruled,vlined,linesnumbered]{algorithm2e}
\let\oldnl\nl
\newcommand{\nonl}{\renewcommand{\nl}{\let\nl\oldnl}}%

\usepackage[utf8]{inputenc}
\usepackage{graphicx}
\usepackage{wrapfig}

\usepackage{amsfonts}
\usepackage{amsmath}
\numberwithin{equation}{section}
\usepackage{mathtools}

\usepackage{makecell}

\usepackage[ruled,vlined,linesnumbered]{algorithm2e}
\let\oldnl\nl

\newcommand{\tenprod}{\circ}
\newcommand{\ten}[1]{\mathcal{#1}}
\newcommand{\mat}[1]{\mathbf{#1}}
\renewcommand{\vec}[1]{\mathbf{#1}}

\newcommand{\PGRAPH}[1]{\noindent\textbf{#1}}

\newcommand{\HAT}[1]{\widetilde{}}

\usepackage{lipsum}

\def\bA{\mathbf{A}}
\def\bU{\mathbf{U}}

\def\bF{\mathbf{F}}

\def\bS{\mathbf{S}}

\def\I{\mathbf{I}}

\newcommand{\vertiii}[1]{{\left\vert\kern-0.25ex\left\vert\kern-0.25ex\left\vert #1 
    \right\vert\kern-0.25ex\right\vert\kern-0.25ex\right\vert}}

\def \diag{{\rm diag}}

\newcommand{\tentt}[1]{\ten{#1}^{TT}}
\SetKwComment{Comment}{/* }{ */}
\makeatletter
\def\ps@pprintTitle{%
 \let\@oddhead\@empty
 \let\@evenhead\@empty
 \def\@oddfoot{}%
 \let\@evenfoot\@oddfoot}
\makeatother

\begin{document}

\begin{frontmatter}

  \title{Tensor Network Space-Time Spectral Collocation Method for Solving the Nonlinear Convection Diffusion Equation}

  \author[TDIV]{Dibyendu Adak}
  \author[TDIV]{M. Engin Danis}
  \author[TDIV]{Duc P. Truong}
  \author[TDIV]{Kim \O. Rasmussen}
  \author[TDIV]{Boian S. Alexandrov}

  \address[TDIV]{Theoretical Division,
    Los Alamos National Laboratory, Los Alamos, NM 87545, USA
  }
  \begin{keyword}
      Spectral Collocation Methods, Chebyshev Polynomial, Exponential Convergence, Tensor Train, TT Rounding, TT Cross Interpolation
      \MSC 65N12 \sep 65N25 \sep 15A23 \sep 15A69
  \end{keyword}

\begin{abstract}
Spectral methods provide highly accurate numerical solutions for partial differential equations, exhibiting exponential convergence with the number of spectral nodes. Traditionally, in addressing time-dependent nonlinear problems, attention has been on low-order finite difference schemes for time discretization and spectral element schemes for spatial variables. However, our recent developments have resulted in the application of spectral methods to both space and time variables, preserving spectral convergence in both domains. Leveraging Tensor Train techniques, our approach tackles the curse of dimensionality inherent in space-time methods. Here, we extend this methodology to the nonlinear time-dependent convection-diffusion equation. Our discretization scheme exhibits a low-rank structure, facilitating translation to tensor-train (TT) format. Nevertheless, controlling the TT-rank across Newton's iterations, needed to deal with the nonlinearity, poses a challenge, leading us to devise the ``Step Truncation TT-Newton" method. We demonstrate the exponential convergence of our methods through various benchmark examples. Importantly, our scheme offers significantly reduced memory requirement compared to the full-grid scheme.
\end{abstract}
\end{frontmatter}


\section{Introduction}

Spectral techniques have been applied for many years to numerically solve  partial differential equations (PDEs). The strength of spectral methods is their exponential convergence with the number of spectral nodes \cite{funaro1997spectral,hussaini1989spectral}. For time-dependent problems, the most well-known schemes combine finite difference methods for the approximation of the temporal variable and spectral methods for the spatial variables. However, these schemes are sub-optimal because the temporal discretization error dominates the spatial discretization error \cite{thomee2007galerkin}. To mitigate this issue, we have previously \cite{adak2024tensor} developed spectral collocation discretization for the linear convection-diffusion-reaction (CDR) equation in space-time format. In \cite{lui2017legendre}, the authors discussed the extension of the space-time spectral collocation method to the nonlinear time-dependent CDR problem. In this case, the discretization of the PDE results in a large system of nonlinear algebraic equations whose numerical solution typically requires Newton's method \cite{knoll2004jacobian}. Again, the nonlinear problem can either be solved using a time-stepping scheme or in space-time format. Time-stepping schemes requires several consecutive iterations of the system of nonlinear equation at each time step, which is very time consuming. On the other hand, the space-time approach results in enormous system of equations, which is very computationally demanding to solve and suffers from the curse of dimensionality \cite{bellman1966dynamic}. The curse of dimensionality manifest itself in the increase of the computational complexity and memory requirements with each additional dimension of the problem, which makes these numerical methods impractical. Importantly, the challenge of the curse of dimensionality remains a significant obstacle in all high-fidelity numerical computations, even in the era of exascale high-performance computing.  

Recently, tensor networks (TNs) \cite{bachmayr2016tensor} have gained attention as a promising method to address the curse of dimensionality. TNs, which generalize tensor factorization \cite{cichocki2014tensor}, mitigate this issue by reorganizing high-dimensional data into networks of low-dimensional small tensors. Originally developed in theoretical physics, TN methods now also show great potential for providing accurate and efficient numerical solutions to high-dimensional large partial differential equations (PDEs). An example of a simple tensor network is the Tensor Train (TT) decomposition \cite{oseledets2011tensor,oseledets2010approximation}. Recent applications of TNs in solving PDEs include efficient and accurate solutions for the Poisson equation \cite{khoromskij2011qtt}, the time-independent Schrödinger equation \cite{gelss2022solving}, the Poisson-Boltzmann equation \cite{benner2021regularization}, Maxwell equations \cite{manzini2023tensor}, Vlasov–Maxwell equations \cite{ye2023quantized}, neutron transport equation \cite{truong2023tensor}, time-dependent convection-diffusion-reaction equation \cite{adak2024tensor}, incompressible fluid dynamics \cite{peddinti2024quantum}, compressible flow Euler equations \cite{danis2024tensortrain}, and others.

However, the extension of the application of TTs to nonlinear PDEs is not straightforward. The primary issues is the growth of the rank of the solution in the Newton method iterations that are required to solve the nonlinear problem at each time point. This issue is caused by the tendency of the Newton method's iteration values to 
leave the low-rank tensor manifolds. This rank growth limits the advantages of tensorization of the PDE. An approach to mitigate this problem is to utilize { \it dynamic low-rank techniques} \cite{dektor2024collocation}. Here, we introduce a new algorithm based on the step-truncation method \cite{rodgers2022adaptive}, which modifies the TT-rank in each Newton iteration, to constrain the solution to the low-rank tensor manifold. 

The outline of this paper is: In Section~\ref{sec: fullgrid}, the governing equations and the full grid discretization are discussed. In Section~\ref{sec: tensor train}, tensor train format and related concepts are introduced. In Section~\ref{sec:tensorization}, we describe in detail the tensorization of the space time spectral collocation scheme, and introduce our {\it step truncation TT-Newton} method. The numerical experiments demonstrating the proposed method are presented in Section~\ref{sec:numerical}.

\section{Model and Numerical Discretization: Full grid scheme}
\label{sec: fullgrid}
\subsection{Nonlinear Convection-Diffusion (NCD) Equation}
\label{eq:model:problem}

We are interested in the following NCD problem:
\begin{align}
    \frac{\partial u}{\partial t}- a(u) \Delta u + \mathbf{b}(u)\cdot \nabla u&= f(u)  \quad \forall (t,\mathbf{x}) \in [0,T] \times \Omega, \label{Model:Eqn1}\\
    u(t,\mathbf{x})&=g(t,\mathbf{x}) \quad \text{on}  \quad [0,T]\times \partial \Omega, \label{Model:Eqn2}\\
    u(0,\mathbf{x})&= h(\mathbf{x}) \quad \text{on} \; \label{Model:Eqn3}\Omega,
\end{align}
where $a(u)$ and $\Vec{b}(u)$ are the nonlinear diffusion and convection coefficients, respectively, while $f(u)$ is the force function.
Furthermore, $\Omega\subset\mathbb{R}^3$ is a three-dimensional (3D), open
parallelepiped domain with boundary $\partial \Omega$; $h(\mathbf{x})$
is the initial condition, (IC), and $g(t,\mathbf{x})$ are the boundary
conditions (BC).
Here, we denote the 3D vectors as, $\Vec{x}=(x,y,z)$, and
matrices also in bold font.

\subsection{Chebyshev Collocation Method}
\label{sec:spacetimeapproach}

Here, we utilize orthogonal Chebyshev polynomials as a global set of basis functions to expand the approximate solution of the NCD equation. Subsequently, we enforce the NCD equation at specific discrete points within the domain, referred to as collocation points \cite{hussaini1989spectral}. The result is a system of nonlinear equations of the unknown expansion coefficients. Based on modified Chebyshev polynomials, $l_i(x)$, (see \cite{funaro1997spectral} for definition), the expansion of approximate solution, $u_h$, is

\begin{equation}
\label{Int_legen}
    u_h(x):= \sum_{j=0}^N u(x_j) l_j(x) \quad \text{where } l_j(x_i) = 
  \left \{
  \begin{aligned}
    & 1 \qquad \text{if } \ i=j, \\
    & 0 \qquad \text{if } \ i \neq j.
  \end{aligned}
  \right.
\end{equation}

The first and second order derivative matrices are constructed using
the matrix representation of the one-dimensional derivative,
$\mathbf{S}^{(N+1) \times (N+1)}_x$, as follows:
\begin{equation}
  \label{derivative_matrix}
  \begin{aligned}
    &\left(\frac{\partial}{\partial x}\right)_{ij} \rightarrow 
  (\mathbf{S}_x)_{i j}:= \frac{d}{dx}l_j(x)\Big|_{x_i},\\
  &\left(\frac{\partial^2}{\partial x^2}\right)_{ij} \rightarrow (\mathbf{S}_{xx})_{ij} = \sum_{s=0}^{N} (\mathbf{S}_x)_{is} (\mathbf{S}_x)_{sj}.
  \end{aligned}
\end{equation}
We refer to \cite{funaro1997spectral} for the
derivation of the last expression.

\if
As a numerical method for PDEs, the Chebyshev collocation method
demonstrates \cite{funaro1997spectral}:
$(i)$ High accuracy, based on its exponential convergence rate with the
number of degrees of freedom;
$(ii)$ Accuracy for capturing complex variations of the solution, based
on the use of smooth basis functions;
$(iii)$ A global approximation, since each basis function affects the
entire domain and is suitable for problems with non-local
dependencies.
Chebyshev collocation method requires regular grids, that is, it works
mainly on Cartesian grids that are tensor products of 1-D domain
partitions, which is suitable to our tensor network formats.
\fi

\subsection{Space Time Matrix Discretization of the NCD equation}
First, we introduce the \emph{space-time} matrix operators: $\frac{\partial}{\partial t}
\rightarrow$ $\mathbf{A}_{t}$; $a(u) \Delta $
$\rightarrow$ $\mathbf{A}_{D}(\mathbf{U})$;  and $\mathbf{b}(u) \cdot \nabla $
$\rightarrow$ $\mathbf{A}_{C}(\mathbf{U})$, and $f(u) \rightarrow \mathbf{F}$ respectively.
Upon employing these matrices, Eq.(\ref{Model:Eqn1}) results in
the following linear system,
\begin{equation}
\boldsymbol{\mathcal{S}}(\bold{U})\bold{U} = \Big (\bold{A}_{t}+ \bold{A}_{D}(\bold{U})+\bold{A}_{C}(\bold{U}) \Big)  \bold{U}=\bold{F}(\bold{U}),
  \label{mat_eq}
\end{equation}
where $\bold{U}$ is a vector corresponding to the
solution, $u_h(t, \bold{x})$, and $\bold{F}(\bold{U})$ is the nonlinear forcing term.
The matrices $\bold{A}_{t}$, $\bold{A}_{D}(\bold{U})$, and $\bold{A}_{C}(\bold{U})$ are of size $(N+1)^4 \times (N+1)^4$ for $N^{th}$ order
spectral collocation method, while $\bold{U}$, $\bold{F}((\bold{U}))$ are column
vectors of size $(N+1)^4 \times 1$, and are designed to incorporate the
boundary conditions in Eq.(\ref{Model:Eqn2}) and Eq (\ref{Model:Eqn3}).

\PGRAPH{$\bullet$ Time Discretization on Chebyshev Grid}: 
To have global exponential convergence, we apply Chebyshev spectral
collocation method for discretization of both spatial and temporal
variables.
To accomplish this, we construct a one-dimensional differential
operator, $\frac{\partial}{\partial t}$ in matrix form, $(\bold{S}_t
)_{ij}$, on a temporal Chebyshev grid, with collocation points $t_0,
t_1\ldots, t_{N}$, as follows (see Eq.\eqref{derivative_matrix}).
\begin{equation}
  \label{time_collocation}
  (\bold{S}_t)_{ij}=\frac{dl_{j}(t)}{dt} \Big |_{t_i}, \quad 0\leq i, j \leq N.
\end{equation}
Hence, the  nonlinear system, in the space-time collocation method is,
\begin{equation}
  \Big (\bold{S}_t \otimes \bold{I}_{\text{space}} +\bold{I}_t \otimes \bold{S}_{\text{space}}(\bold{U}) \Big)\bold{U} = \bold{F}(\bold{U}), 
  \label{cheb_lin_syst}
\end{equation}
where (see Eq. \ref{mat_eq}) $\mat{A}_t=\bold{S}_t \otimes \bold{I}_{\text{space}}$, and $\otimes$ is the Kronecker product.
To construct the space-time operator of the NCD equation on the
Chebyshev grid we also need to construct the nonlinear space operator,
$\bold{S}_{\text{space}}(\bold{U})$. $\bold{I}_{\text{space}}$ and $\bold{I}_t$ are identity matrices with the dimensions of the space and time variables, respectively. For simplicity we will write the discretization of the nonlinear system as,
\begin{equation}
    \bold{G}(\bold{U})=\bold{0}
    \label{non_lin_sys}
\end{equation}
where, 
\begin{equation*}
    \bold{G}(\bold{U}):= \Big (\bold{S}_t \otimes \bold{I}_{\text{space}} +\bold{I}_t \otimes \bold{S}_{\text{space}}(\bold{U}) \Big)\bold{U} - \bold{F}(\bold{U})
\end{equation*}

\PGRAPH{$\bullet$ Discretization of the Nonlinear Diffusion Operator:} The Laplacian, $\textbf{L}$, on the spatial Chebyshev grid, is constructed as follows,
\begin{equation}
  \begin{split}
    \label{Spectral_Diff}
    \bold{L} =  \bold{I}_{t}\otimes \bold{S}_{xx} \otimes \bold{I}_{y}\otimes  \bold{I}_{z}+\bold{I}_{t} \otimes \bold{I}_{x} \otimes \bold{S}_{yy}  \otimes\bold{I}_{z}  + \bold{I}_{t} \otimes \bold{I}_{x} \otimes  \bold{I}_{y} \otimes\bold{S}_{zz}.
  \end{split}
\end{equation} 
The nonlinear operator $\bA_D (\bold{U})$ is the product of the nonlinear diffusion coefficient, $a(u)$, and the Laplacian, $\bold{L}$:
\begin{equation}
     \bA_D (\bold{U})= \diag(a(\bold{U})) \bold{L},
\end{equation}
\PGRAPH{$\bullet$ Discretization of the Convection Term}:
Here, we focus on the matricization of the nonlinear convection term,
$\bold{A}_C (\bold{U})$, with the
convective function, $\bold{b}(u)$, which we assume in the form
\begin{equation}
  \vec{b}(u)=[b^x(u)~b^y(u)~b^z(u)].
\end{equation}
Then, we construct the nonlinear convection term $\bold{A}_C(\bold{U})$
\begin{equation}
  \label{convec:mat}
  \begin{split}
    \bold{A_C(U)}= \ & \diag(\bold{B}^x(\bold{U})) \left(\bold{I}_{t} \otimes \bold{S}_{x} \otimes \bold{I}_{y}\otimes\bold{I}_{z}\right) + \ldots\\
    & \diag(\bold{B}^y(\bold{U}))\left(\bold{I}_{t} \otimes \bold{I}_{x}\otimes \bold{S}_{y}  \otimes\bold{I}_{z}\right) + \ldots \\
    & \diag(\bold{B}^z(\bold{U}))\left(\bold{I}_{t} \otimes \bold{I}_{x} \otimes  \bold{I}_{y}\otimes \bold{S}_{z}\right),
  \end{split}
\end{equation}
where $\bold{B}^x,\bold{B}^y$, and $\bold{B}^z$ are vectors of size
$(N+1)^4$ containing the evaluation of the functions $b^x(u),\ b^y(u)$ and
$b^z(u)$ on the Chebyshev space-time grid.

\subsection{Initial and Boundary Conditions on Space-Time Chebyshev Grids} 
\label{Bd_cond}

\begin{figure}[!t]
  \centering
  \includegraphics[width=0.5\textwidth]{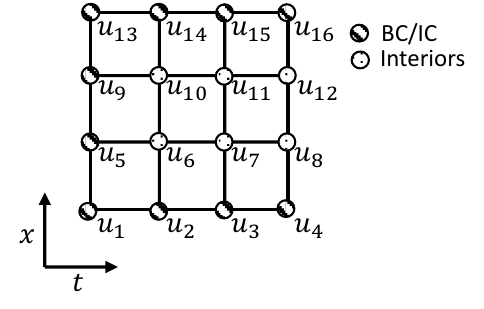}
  \caption{1D Space time grid with $N=4$ collocation nodes.}
  \label{fig:ST-2d-grid}
\end{figure}

So far, we have neglected the boundary conditions (BC) and the
initial condition (IC), given in \eqref{Model:Eqn3}.
In the space-time method, we consider the IC equivalent to the BC.
The nodes of the Chebyshev grid are split into two parts: $(i)$ BC/IC
nodes, and $(ii)$ interior nodes, where the solution is unknown.
Let $i^{\text{Bd}}$, and $i^{\text{Int}}$ be the set of multi-indices,
introduced in \cite{dolgov2021guaranteed},
for the BC/IC and interior nodes, respectively.

We impose the boundary and initial conditions explicitly by enforcing
$u_h(t,x)=g(t,x)$ on the BC or $u_h(0,x)=h(x)$ on the IC. Using this, we reduce the nonlinear system \eqref{mat_eq} for all nodes into a smaller
system representing only the interior nodes.
To make it clear, we consider below a simple example with $N = 4$
collocation points in two dimensions $(t,x)$.
The set of Chebyshev nodes can be denoted by multi-indices as,
$\bold{U}:=\{u_1,u_2,\ldots,u_{16} \} $ (Figure
\ref{fig:ST-2d-grid}). After imposing the BC/IC conditions, the nonlinear
system becomes, 

\begin{equation}
  \label{Full_system}
  \begin{split}
    &u_1 = g (t_{0},{x}_{0}),\\
    &\quad \vdots \\
    &u_4 = g (t_{3},{x}_{0}),\\
    &u_5 = h(x_1), \\
    &s_{6,1}(u_6) u_1 + s_{6,2}(u_6) u_2+ \ldots +s_{6,15}(u_6)u_{15} + s_{6,16}(u_6)u_{16} -f(u_6) =0,\\
    &\quad \vdots \\
    &s_{8,1}(u_8) u_1 + s_{8,2}(u_8) u_2+ \ldots +s_{8,15}(u_8)u_{15} + s_{8,16}(u_8)u_{16} - f(u_8)=0,\\
    &u_9 = h(x_{2})\\
    &s_{10,1}(u_{10}) u_1 + s_{10,2}(u_{10}) u_2+ \ldots +s_{10,15}(u_{10})u_{15} + s_{10,16}(u_{10})u_{16} -f(u_{10})=0,\\
    &\quad \vdots \\
    &s_{12,1}(u_{12}) u_1 +s_{12,2}(u_{12}) u_2+ \ldots +s_{12,15}(u_{12})u_{15} + s_{12,16}u_{16} - f(u_{12})  = 0,\\
    & \quad \vdots \\
    &u_{13} = g (t_{0},{x}_{3}),\\
    &\quad \vdots \\
    &u_{16} = g (t_{3},{x}_{3}).\\
  \end{split}
\end{equation}
Rearranging these equations, we find that the unknown values $(u_6,u_7,u_8,u_{10},u_{11},u_{12})$ associated with the
six interior nodes satisfy the following system of equations
\begin{multline}
  \qquad
  s_{l,6}(u_l) u_6+s_{l,7}(u_l)u_7+s_{l,8}(u_l)u_{8} + s_{l,10}(u_l)u_{10} +s_{l,11}(u_l)u_{11} +s_{l,12}(u_l)u_{12}  \\
  -f(u_l)  + s_{l,1}(u_l)g(t_0,x_0)+s_{l,2}(u_l)g(t_1,x_0)+\ldots = 0
  \label{eqn:sub_system_BC}
\end{multline}
with $l \in \{6,7,8,10,11,12\}$. 
Generalizing to three spatial dimensions and using matrix form, Eq. (\ref{eqn:sub_system_BC}) becomes
\begin{equation}
    \bold{S}^{reduced}(\bold{U}_{int})\bold{U}_{int}-\bold{F}(\bold{U}_{int})+\bold{F}^{BC}(\bold{U}_{int}) = 0,
    \label{red_non_lin_sys}
\end{equation}
where,
\begin{equation*}
\bold{S}^{reduced}(\bold{U}_{int}):= \bold{A}^{reduced}_{t}+ \bold{A}^{reduced}_{D}(\bold{U}_{int})+\bold{A}^{reduced}_{C}(\bold{U}_{int}),
\end{equation*}
and $\bF^{BC}(\bU_{int})$ is the boundary terms that correctly represent boundary and initial conditions in the nonlinear system for the interior nodes. Then the nonlinear system (\ref{non_lin_sys}) then becomes,
\begin{equation}
\textbf{G}^{reduced}(\bold{U}_{int}) = 0.
\label{red_res}
\end{equation}

\subsection{Newton's Method with Line Search}
\label{Newton:Method}
Given that the reduced system in Eq.\eqref{red_non_lin_sys} is a nonlinear system, and that its Jacobian can be computed analytically, we choose to use Newton's method to approximate its solution.
Newton's method is an iterative algorithm for finding the roots of a the residual function the nonlinear system (\ref{red_res}). The iteration derives from a Taylor expansion around the current iteration $k$ of the solution, $\bold{U}_{int}^{k}$,
\begin{equation}
    \bold{G}^{reduced}(\bold{U}_{int}^{k+1}) = \bold{G}^{reduced}(\bold{U}_{int}^k) + \textbf{J}(\bold{U}_{int}^{k}) (\bold{U}_{int}^{k+1} - \bold{U}_{int}^k) + \text{higher order terms}. 
\end{equation}
Here, $\bold{J}$ is the Jacobian matrix corresponding to the nonlinear system $\textbf{G}^{reduced}(\bold{U}_{int})$, and $k$ is the iteration index. 
If the higher order terms is neglected and the right-hand side is set to zero, we arrive at Newton iteration for a given initial guess, $\bU_0$ as a linear system relative to the correction, $\boldsymbol{\delta}^k$, leading to the next step, $\bold{U}_{int}^{k+1}$,
\begin{equation}
\begin{aligned}
    \bold{J}(\bold{U}_{int}^k) \boldsymbol{\delta}^k &= - \bold{G}^{reduced}(\bold{U}_{int}^k) \\
    \bold{U}_{int}^{k+1} &= \bold{U}_{int}^{k} + \boldsymbol{\delta}^k.
\end{aligned}
\label{newton_it}
\end{equation}
We solve this linear system numerically.
To improve the convergence robustness of Newton's method, we employ a standard globalization technique, called line search method \cite{kelley1995iterative}. In the line search method, a factor ${s}$ is used to decide how far the update should be in the direction of $\boldsymbol{\delta}^k$, and  $\bU_{int}^{k+1} = \bU_{int}^k + s\boldsymbol{\delta}^k$.

\if

\section{Tensor Train Representation}
\label{sec: tensor train}
In this section, we describe the tensorization process for the space-time collocation discretization of the NCD equations. The tensor network we use here is tensor train (TT) format. We also utilize the cross interpolation method to approximate the TT format of coefficients and input functions. For a more comprehensive understanding of notation and concepts, we refer the reader to the following references:
~\cite{adak2024tensor,kolda2009tensor,oseledets2011tensor,dolgov2021guaranteed}.
We also propose a Newton method in TT format for solving nonlinear tensor equations, utilizing the step-truncation method described in ~\cite{rodgers2023implicit,rodgers2022adaptive}. The algorithm involves projecting the updated solution onto a low-rank tensor manifold through a rounding procedure at each Newton iteration, maintaining the low-rank structure, and thereby, enhancing computational efficiency.

\subsection{Tensor Train}
Tensor train network, or TT-format of a tensor, \cite{oseledets2010tt} represents a tensor as a product of cores, which are either matrices or three-dimensional tensors. Since the tensors in our numerical schemes are four-dimensional (4D), we focus exclusively on the TT format for 4D tensors. Specifically, the TT approximation $\mathcal{X}^{TT}$ of a 4D tensor $\mathcal{X}$ is defined as follows:
\begin{equation}
    \label{eq:TT_def}
    \ten{X}_{TT}(i,j,k,l) = \sum_{\alpha_1=1}^{r_1}\sum_{\alpha_2=1}^{r_2}\sum_{\alpha_3=1}^{r_3}\ten{G}_1(1,i,\alpha_1)\ten{G}_2(\alpha_1,j,\alpha_2)\ten{G}_3(\alpha_2,k,\alpha_3)\ten{G}_4(\alpha_3,k,1) + \varepsilon(i,j,k,l),
\end{equation}
where the error, $\varepsilon$, is a 4D tensor with the same dimensions as $\ten{X}$. The integers $[r_1, r_2, r_3]$ are known as TT-ranks. Each TT-core, $\ten{G}k$, depends on only one index of $\ten{X}_{TT}$, indicating that the TT format represents a discrete separation of variables.

Equivalently, the elements of $\ten{X}^{TT}$ can be represented as a product of vectors and matrices as follows:
\begin{equation}
    \label{eq:TT_def_2}
    \ten{X}_{TT}(i,j,k) = G_1(i)G_2(j)G_3(k)G_4(l) + \varepsilon(i,j,k),
\end{equation}
where $G_1(i)$ is a vector of size $1 \times r_1$, $G_2(j)$, and $G_3(k)$ are matrices of size $r_1 \times r_2$, and $r_2 \times r_3$ respectively. $G_4(k)$ is a vector of size $r_3 \times 1$.
\begin{figure}
    \centering
    \includegraphics[width = 0.8\textwidth]{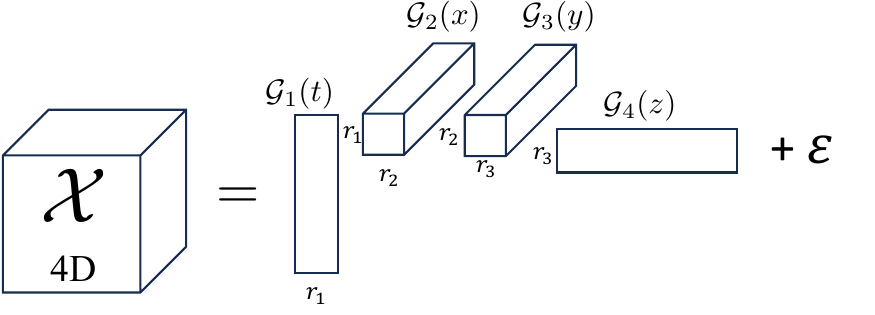}
    \caption{TT format of a 4D tensor $\ten{X}$, with TT cores $\ten{G}_1,\ \ten{G}_2,\ \ten{G}_3$, TT-ranks $\mathbf{r} = (r_1,r_2)$, and approximation error presented in the 4D tensor, $\varepsilon$.}
    \label{fig:TT-format}
\end{figure}
\subsection{TT Rounding}
\label{rounding}
When a tensor is already in TT format, denoted as $\ten{X}_{TT}$ with TT-ranks $r_k$, the TT rounding procedure is employed to determine a more compact TT representation, $\ten{Y}_{TT}$, with TT-ranks ${r'}_k \leq r_k$, while ensuring a prescribed accuracy $\varepsilon_{TT}$, such that:
\begin{equation}
    \Vert \ten{X}_{TT} - \ten{Y}_{TT}\Vert \leq \varepsilon_{TT}\Vert \ten{X}_{TT} \Vert.
\end{equation}

This process is commonly known as TT rounding, truncation, or recompression. The TT rounding algorithm builds on the relationship between truncation errors of the TT format and its cores \cite[Theorem 2.2]{oseledets2010tt}, efficiently employing QR and truncated singular value decompositions (SVD) at the TT core level \cite{oseledets2011tensor}. The TT rounding algorithm consists of two passes along the cores. During the first pass, the tensor is orthogonalized either from left to right or right to left, followed by a sweep in the opposite direction during the second pass. The TT-ranks are successively reduced by truncating the SVD of the matricized cores.

\subsection{TT Cross Interpolation} 
In the context of this paper, we need a method to construct TT format of coefficient functions, initial conditions, boundary conditions and source function. A straightforward way to do this is using the TT SVD algorithm, which relies on a series of singular value decompositions (SVD) carried out on the unfolding matrix of a tensor \cite{oseledets2011tensor}. While known for its efficiency, this algorithm requires access to the entire tensor, which can be impractical or even impossible for large tensors. To address this challenge, the cross interpolation algorithm, also known as TT-cross, was introduced in \cite{oseledets2010tt}.
Rather than employing SVD, TT-cross utilizes an approximate variant of the skeleton/CUR decomposition \cite{mahoney2009cur}. The CUR decomposition approximates a matrix $A \approx CUR$, where C and R represent a selection of columns and rows from A, respectively, and U is the inverse of the intersection sub-matrix. To identify an effective set of rows and columns, the TT-cross algorithm takes advantage of the MaxVol algorithm, which is derived from the Maximum Volume Principle \cite{goreinov2010find}.
CUR decomposition does not require access to the complete tensor; instead, it only requires a function or routine to compute tensor elements on-the-fly.
Nevertheless, directly extending CUR to high-dimensional tensors remains computationally expensive. Consequently, alternative heuristic tensor network optimization-based algorithms for cross interpolation have been developed, including the Alternating Linear Scheme \cite{holtz2012alternating}, Density Matrix Renormalization Group (DMRG) \cite{savostyanov2011fast}, and Alternating Minimal Energy (AMEn) \cite{dolgov2014alternating}.
These iterative algorithms aim to approximate the TT format of a large tensor $\ten{A}$ with a specified accuracy $\varepsilon_{cross}$, such that:
\begin{equation}
    \Vert \ten{A} - \ten{A}_{TT} \Vert \leq \varepsilon_{cross} \Vert \ten{A} \Vert.
\end{equation}

\section{Tensorization of NCD}
\label{sec:tensorization}
Omitting for simplicity the subscripts $\emph{int}$ and superscripts $\emph{reduced}$, Eq.~\eqref{red_non_lin_sys} for all interior nodes becomes
\begin{equation}
\label{eq:last_reduced}
\bold{A}_{t}\bU + \bold{A}_{D}(\bU)\bU + \bold{A}_C(\bU)\bU - \bold{F}(\bold{U}) + \bF^{BC}(\bU) = 0.\\
\end{equation}

The tensorization is a process of representing Eq.~\eqref{eq:last_reduced} in TT format, by constructing the TT-format of the discrete operators, for more details please see Ref. \cite{truong2023tensor}. In TT format, the nonlinear system is,
\begin{equation}
\label{eq:TT-non-sys}
\tentt{A}_{t}\tentt{U} + \tentt{A}_D(\tentt{U})\tentt{U} + \tentt{A}_{C}(\tentt{U})\tentt{U} - \ten{F}(\tentt{U}) + \ten{F}^{bc}(\tentt{U}) = 0.
\end{equation}



Given that time, diffusion and convection operators have Kronecker product structure, their TT format can be exactly formed by using the component matrices as the TT cores. To construct the TT format of the operators acting only on the interior nodes, some index sets need to be defined:
\begin{equation*}
  \begin{split}
    & \ten{I}_t = 1:N \ \text{index set for time variable},\\
    & \ten{I}_s = 1:(N-1) \ \text{index set for space variable}.
  \end{split}
\end{equation*}

\PGRAPH{$\bullet$ TT-matrix Time Operator}, $\ten{A}^{TT}_t$: From the Kronecker product of the time operator in~\eqref{cheb_lin_syst}, its TT-matrix format acting on the interior nodes is constructed as: 

\begin{equation}
  \ten{A}^{TT}_t = \bold{S}_t(\ten{I}_t,\ten{I}_t)\tenprod \bold{I}_{N-1} \tenprod \bold{I}_{N-1} \tenprod \bold{I}_{N-1},
  \label{eqn:At_TT_op}
\end{equation}
where $\bold{I}_{N-1}$ is the identity matrix of size
$(N-1)\times(N-1)$, and $\tenprod$ is the tensor product operator.
  
\PGRAPH{$\bullet$ TT-matrix of the Nonlinear Diffusion Operator}, $\tentt{A}_D(\tentt{U})$: 
From its structure in Eqn.~\ref{Spectral_Diff}, the Laplace operator in TT-matrix format is constructed as:
\begin{equation}
  \begin{split}
      \tentt{L}=\bold {I}_{N}\tenprod \bold{S}_{xx}(\ten{I}_s,\ten{I}_s) \tenprod \bold{I}_{N-1}\tenprod \bold{I}_{N-1}
      + &\bold{I}_{N} \tenprod \bold{I}_{N-1} \tenprod \bold{S}_{yy}(\ten{I}_s,\ten{I}_s)  \tenprod\bold{I}_{N-1}\\
      + &\bold{I}_N \tenprod \bold{I}_{N-1} \tenprod  \bold{I}_{N-1} \tenprod\bold{S}_{zz}(\ten{I}_s,\ten{I}_s).
  \end{split}
\end{equation}
Then the coefficient function $a(u)$ is incorporated as, 
\begin{equation}
\tentt{A}_D(\tentt{U}) = \diag(a(\tentt{U}))\tentt{L}.
\end{equation}
The $\diag()$ operation for TT is described in~\ref{app:func_eval}.

\PGRAPH{$\bullet$ TT-matrix of the Nonlinear Convection Operator}, $\tentt{A}_C(\tentt{U})$:
From its structure in Eqn.(\ref{convec:mat}), the convection operator in TT-matrix format is constructed as:
\begin{equation}
  \begin{split}
    \nabla^{TT}_x = \bold {I}^{N}\tenprod \bold{S}_{x}(\ten{I}_s,\ten{I}_s) \tenprod \bold{I}^{N-1}\tenprod \bold{I}^{N-1},\\[0.5em]
    \nabla^{TT}_y = \bold {I}^{N}\tenprod \bold{I}^{N-1}\tenprod\bold{S}_{y}(\ten{I}_s,\ten{I}_s) \tenprod \bold{I}^{N-1}, \\[0.5em]
    \nabla^{TT}_z = \bold {I}^{N} \tenprod \bold{I}^{N-1}\tenprod \bold{I}^{N-1} \tenprod\bold{S}_{z}(\ten{I}_s,\ten{I}_s).\\
  \end{split}
\end{equation}
Then the nonlinear convection operator $\tentt{A}_C(\tentt{U})$ becomes 
\begin{equation}
  \tentt{A}_C(\tentt{U}) =
  \diag(b^x(\tentt{U}))\nabla_x^{TT} +
  \diag(b^y(\tentt{U}))\nabla_y^{TT} +
  \diag(b^z(\tentt{U}))\nabla_z^{TT}.
\end{equation}
The construction of the coefficient functions: $\diag(a(\tentt{U}))$,  $\diag(b^x(\tentt{U}))$, loading tensor, $\ten{F}(\tentt{U})$, and the boundary term, $\ten{F}^{BC}(\tentt{U})$ in TT formats are described in the~\ref{sec:app}.
At this stage, we completed the construction of the TT format of the terms in the nonlinear system Eqn.~\eqref{eq:TT-non-sys}, next we will introduce the TT-Newton method to solve TT nonlinear equations.

\subsection{TT-Newton Method With Step Truncation }
\label{TT-Newton}
Our tensorized Newton's method is described below in Algorithm 1. Initial numerical experiments showed that when using TT rounding accuracy with a fixed $\varepsilon_{TT}$ (see Sect \ref{rounding}), the solution rank can increase uncontrollably. This makes the TT-Newton method with a fixed accuracy unusable because the iterations can move the iteration values away from the low-rank tensor manifold.

To maintain low-rank structures throughout the Newton iterations, starting from a relatively large $\varepsilon^0_{TT}$, we allow the TT rounding accuracy in the TT-Newton method, $\varepsilon^k_{TT}$ to vary after each iteration, $k$. A similar approach was introduced in \cite{rodgers2023implicit}, and it was shown that it leads to results similar to those obtained using Dynamic Low Rank Approximation (DLRA) method introduced for solving nonlinear equations in TT format, \cite{lubich2014projector,ceruti2022rank,dektor2024collocation}.
In our implementation we used  two stopping criteria: (i) the norm of nonlinear residual $\text{tol}_\text{res} > \Vert\bold{G}(\ten{U}^{TT,k}) \Vert/ \Vert\bold{G}(\ten{U}^{TT,0}) \Vert$, and (ii) the norm of the update difference $\text{tol}_\text{update} > \left\Vert \dfrac{\boldsymbol{\delta}^{TT,k}}{\ten{U}^{TT,k}}\right\Vert$, which are standard criteria for the Newton's method~\cite{knoll2004jacobian}. Then $\varepsilon^{k}_{TT}$, the TT truncation tolerance of the $k$'th iteration, is updated as,
\begin{equation}
    \label{eq:epsk}
    \varepsilon^k_{TT} = \text{min}\left( \varepsilon^{k-1}_{TT}, \left\Vert\dfrac{\bold{G}(\ten{U}^{TT,k})}{\bold{G}(\ten{U}^{TT,0})}\right\Vert,\left\Vert \dfrac{\boldsymbol{\delta}^{TT,k}}{\ten{U}^{TT,k}}\right\Vert  \right).
\end{equation}

\begin{algorithm}
\caption{Step Truncation TT-Newton Method}\label{alg:two}
\KwData{TT function $\bold{G}(\tentt{U})$, $\bold{J}(\tentt{U})$, initial guess $\ten{U}^{TT,0}$, TT truncation tolerance $\varepsilon^0_{TT}$, and maximum number of iteration $N_{it}$.}
\KwResult{Approximate root $\tentt{U}$ of $\bold{G}(\tentt{U}) = 0$}
\For{k = 0:$N_{it}$}{
Compute and recompress $\bold{G}(\ten{U}^{TT,k})$ and $\bold{J}(\ten{U}^{TT,k})$ with $\varepsilon^{k}_{TT}$.\\
Solve for $\bold{\delta}^{TT,k}$: $\bold{J}(\ten{U}^{TT,k}) \bold{\delta}^{TT,k} = - \bold{G}(\ten{U}^{TT,k})$.\\
\For{$s = 1,0.5,0.25,...$}{
$\ten{U}^{TT,k+1} = \ten{U}^{TT,k} + s\bold{\delta}^{TT,k}$\\
Recompress $\ten{U}^{TT,k+1}$ with $\varepsilon^k_{TT}$. \Comment{Truncation Step}
\uIf{$\bold{G}(\ten{U}^{TT,k+1}) \leq \bold{G}(\ten{U}^{TT,k})$}
{\textbf{break}}}
Check for stopping criteria.\\
Update $\varepsilon^{k+1}_{TT}$ using Eqn.~\eqref{eq:epsk}.
}
\end{algorithm}

\section{Numerical Experiments}
\label{sec:numerical}
Here we conduct several numerical experiments to assess the computational and memory efficiency of the TT space time solver and the step truncation TT-Newton method. In the first experiment, we show the comparison between full grid Newton, step truncation TT-Newton and a fixed-$\varepsilon_{TT}$ TT-Newton solvers.
%
%
In the second experiment, we compare the performance of a full grid solver with TT solvers for solving the NCD equation with a manufactured solution. 
In the third experiment, we show how full grid and TT solvers perform solving a 3D viscous Burger's equation. All our implementations utilize the MATLAB TT-Toolbox~\cite{oseledets2011matlab} on a Linux machine with Xeon Gold 6148 processor.

\subsection{Experiment 1: Step Truncation TT-Newton Method}
In this numerical experiment, we illustrate the need for adaptive truncation tolerance in our TT-Newton method (see, Section \ref{TT-Newton}). We show that allowing the TT-truncation tolerance $\varepsilon_{TT}$ to adaptively vary with the stopping criteria gives us the necessary control over the TT-ranks of the solution. To demonstrate this, we consider the following task: Use the TT-Newton method, to find a 4D tensor, $\ten{Y}$, such that, $q(\ten{Y}) = \exp({-\ten{Y}}) - \ten{Y}^3 - \ten{G}=0$. Here we chose $\ten{G}={\exp(-\ten{Y}_{exact}}) - \ten{Y}_{exact}^3$, making a randomly generated 4D TT-tensor, $\ten{Y}_{exact}$, the exact solution of the problem. 

The equation, $q(\ten{Y}) = 0$, is solved in three ways with the same initial guess and the same stopping criterion $10^{-6}$. The first way is to apply the full grid Newton method, the second way is to use TT-Newton method with a fixed $\varepsilon_{TT} = 10^{-8}$, and the third is to allow $\varepsilon_{TT}$ to vary with the stopping criteria (Eqn.~\eqref{eq:epsk}), starting from, $\varepsilon_{TT,0} = 10^{-1}$.
The full grid Newton method converges after six iterations with the accuracy $5.7\times10^{-8}$. The TT-Newton method with a fixed $\varepsilon_{TT}$ converges after seven iterations, with accuracy $1.07\times10^{-6}$, and the step-truncated TT-Newton method converges after six iterations with accuracy $1.61\times10^{-6}$. Below, we use the compression ratio of a TT $\ten{Y}$ defined as:
\begin{equation*}
    CR(\tentt{Y}) = \dfrac{\text{\# elements in TT}}{\text{\# elements in full tensor}}
\end{equation*}
The elapsed time in seconds and the compression ratio of the TT format for each Newton iteration are given in Figure~\ref{fig:test-1}. From the left panel, it is seen that the full grid Newton method, as expected, requires approximately the same amount of time for each iteration.
In contrast, the TT-Newton method, with a fixed $\varepsilon_{TT}$, requires increasingly more time compared to the full grid Newton method as the iterations progress. This is due to the rapid growth of TT-ranks in each iteration. The right panel illustrates this TT-rank growth, showing that the compression ratio of the TT-solution in the TT-Newton method almost reaches 1, at which point the TT-format provides no compression advantage over full grid tensors. 
Importantly, for the step-truncation TT-Newton method, with adaptive $\varepsilon_{TT}$, the elapsed time is shorter across all iterations compared to the full grid Newton method. This increased efficiency is supported by a favorable compression ratio, shown in the right panel. The results from this test case demonstrate that the step-truncation TT-Newton method successfully maintains the low-rank structures of the solution throughout Newton iterations, and offers a better efficiency than the full grid algorithm.
\begin{figure}
    \centering
    \includegraphics[width=0.9\textwidth]{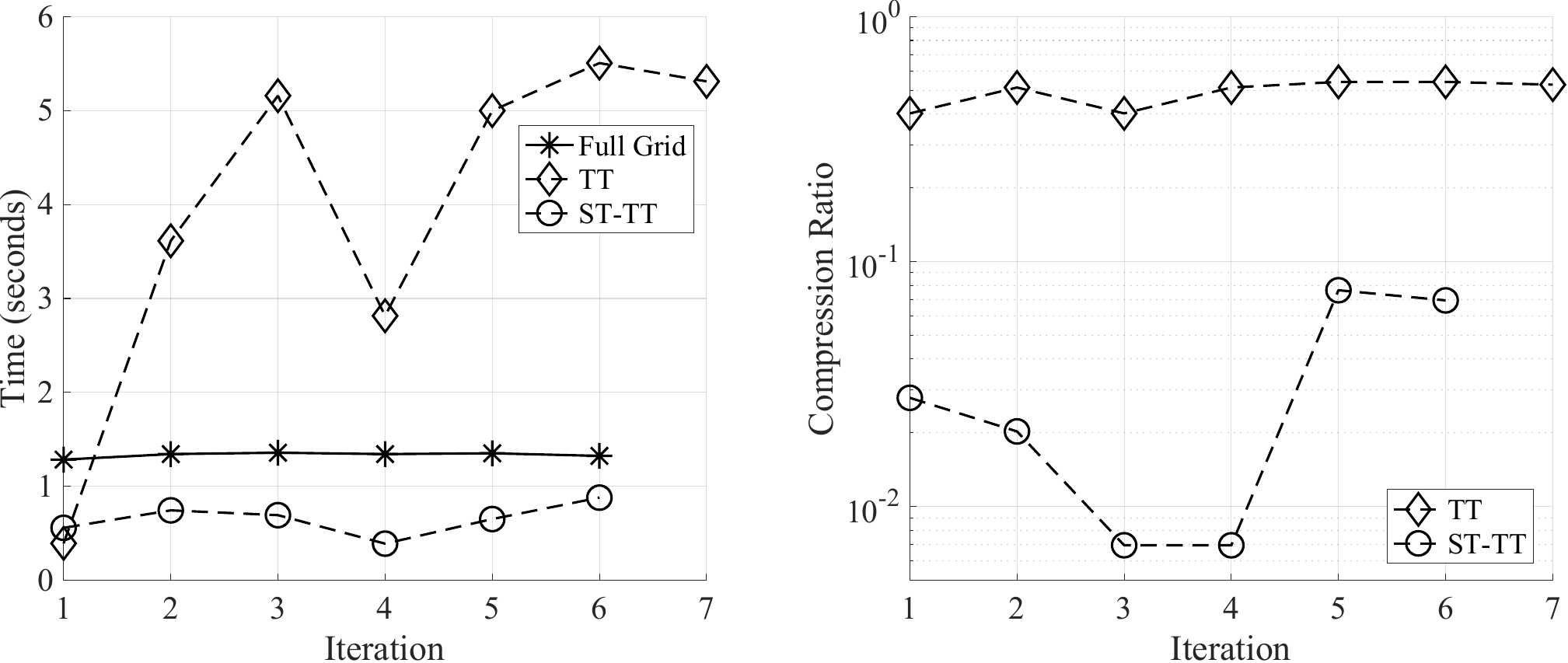}
    \caption{Performance of full grid Newton, TT-Newton with fixed $\varepsilon_{TT}$, and step truncation TT-Newton (ST-TT) solvers across Newton iterations. \textit{Left panel}- Elapsed time for each iteration. ST-TT is most efficient. \textit{Right panel} - Compression ratio of TT and ST-TT, shows that ST-TT successfully maintains the low-rank structures of the solution throughout Newton iterations.}
    \label{fig:test-1}
\end{figure}

\subsection{Experiment 2: Manufactured Solution}
In this experiment, we solve the 3D NCD model problem, Eqn.~\eqref{eq:model:problem}, for a manufactured solution, $u_{exact}(t,\textbf{x}) = \exp({-0.1t})\sin(\pi x)\sin(\pi y) \sin(\pi z)$, and coefficient functions, $a(u) = 1+u^2$, $\mathbf{b}(u) = [u, 1, 1]$, and $f(u) = u - u^3$. The computational domain is, $[0,T]\times \Omega: [0,1]\times[-2,2]^3 $. The inhomogeneous initial and boundary conditions are determined from the manufactured solution. Equation~\eqref{eq:model:problem} is solved in three ways: 

\begin{itemize}
    \item First, using the full grid scheme (see Sec.~\ref{sec: fullgrid}) with Newton stopping criterion of $10^{-6}$. 
    \item Second, using the tensorization scheme (see, Sec.~\ref{sec:tensorization}) with TT rounding accuracy, $\varepsilon_{TT} = 10^{-5}$, and step-truncation TT-Newton method with  stopping criterion of $10^{-6}$.
    \item Finally, using the tensorization scheme with TT-rounding accuracy, $\varepsilon_{TT} = 10^{-9}$, and step-truncation TT-Newton method with stopping criterion of $10^{-8}$.
\end{itemize}
Figure~\ref{fig:test-2} shows the performance of all three cases, including the error convergence plot (left panel), and elapsed time (right panel). The full grid algorithm can only be run for up to 16 collocation points per dimension due to memory limitations. Tensorization of the model problem, with added step-truncation TT-Newton algorithms, ran up to 24 collocation points per dimension for both used accuracies. 
\begin{itemize}
    \item In the left panel, all three experiments show the exponential convergence as expected of the space-time spectral collocation method. The error curves of the TT algorithms plateau as they approach their respective TT rounding accuracy: $\varepsilon_{TT} = 10^{-5}$ and $\varepsilon_{TT} = 10^{-9}$, suggesting that the error of the TT solutions are controlled by $\varepsilon_{TT}$.
   \item In the right panel, the result suggests that when $\varepsilon_{TT}$ is much smaller than the numerical error of the solution, the TT algorithms are inefficient, as the TT algorithms take more time than the full grid simulations. As recently suggested in \cite{danis2024tensortrain}, this likely due to an interaction between the truncation error of the underlying numerical discretization and $\varepsilon_{TT}$. An unnecessarily small $\varepsilon_{TT}$, i.e. $\varepsilon_{TT}$ is smaller than the numerical error, is thought to cause TT recompression to attempt to incorporate the numerical noise into the TT format, which might result in increased TT ranks and slower run times for the TT solver. For a more detailed discussion, we refer the interested readers to Ref.~\cite{danis2024tensortrain}. For the TT algorithm with $\varepsilon_{TT} = 10^{-5}$, this occurs for 8 collocation points; and for the TT algorithm with $\varepsilon_{TT} = 10^{-9}$, it occurs for 8 and 12 collocation points. 
\end{itemize}

\begin{figure}
    \centering
    \includegraphics[width=1.0\textwidth]{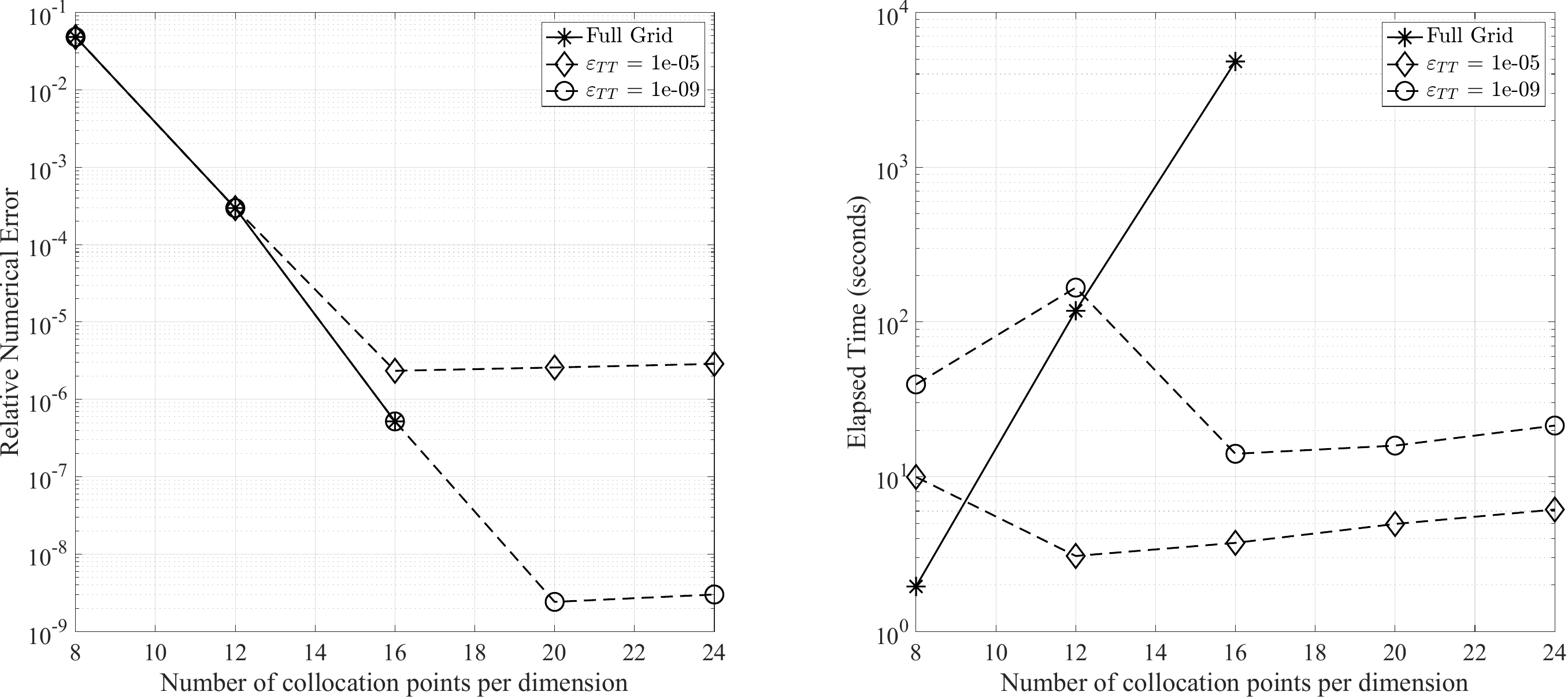}
    \caption{Experiment 2: Left Panel - Relative numerical error versus number of collocation points per dimension. The error curves show the exponential convergence of three schemes. Right Panel - Elapsed time in seconds. TT schemes is more efficient than the full grid scheme.}
    \label{fig:test-2}
\end{figure}

\subsection{Experiment 3: 3D Viscous Burgers' Equation}
Here we consider the 3D viscous Burgers' equation
\begin{equation}
    \dfrac{\partial u}{\partial t} + u(u_x + u_y + u_z) = \Delta u,
\end{equation}
which is equivalent to Eq. (\ref{Model:Eqn1}) for $a(u)=1$, $\mathbf{b}(u)=[u~u~u]$, and $f(u)=0$.
The Burgers' equation resembles the incompressible Navier-Stokes equations. Both equations have nonlinear convective terms and linear viscous terms. However, the Burgers' equation is much simpler. It is a scalar equation and has no coupling with other equations. It also does not have a pressure gradient term.

Following the 1D exact solution presented in \cite{wood2006}, we derive the following exact solution to the 3D viscous Burger's equations:
\begin{equation}
    u(t,x,y,z)=\frac{(2/3)\pi\exp{\left(-\pi^2 t/3\right)}\sin{\left(\pi (x+y+z)/3\right)}}{5+\exp{\left(-\pi^2 t/3\right)}\cos{\left(\pi (x+y+z)/3\right)}}.
\end{equation}
The initial condition and boundary conditions are obtained from the exact solution. The computational domain is chosen as, $[0,T]\times \Omega: [0,1]\times[0,6]^3$. Similar to Experiment 2, the problem is solved with three schemes: full grid scheme with Newton stopping criterion of $10^{-6}$, the TT scheme with $\varepsilon_{TT} = 10^{-5}$ and step-truncation TT-Newton method with stopping criterion $10^{-6}$, and the TT scheme with $\varepsilon_{TT} = 10^{-8}$ and stopping criterion $10^{-7}$.
\begin{figure}
    \centering
    \includegraphics[width=1.0\textwidth]{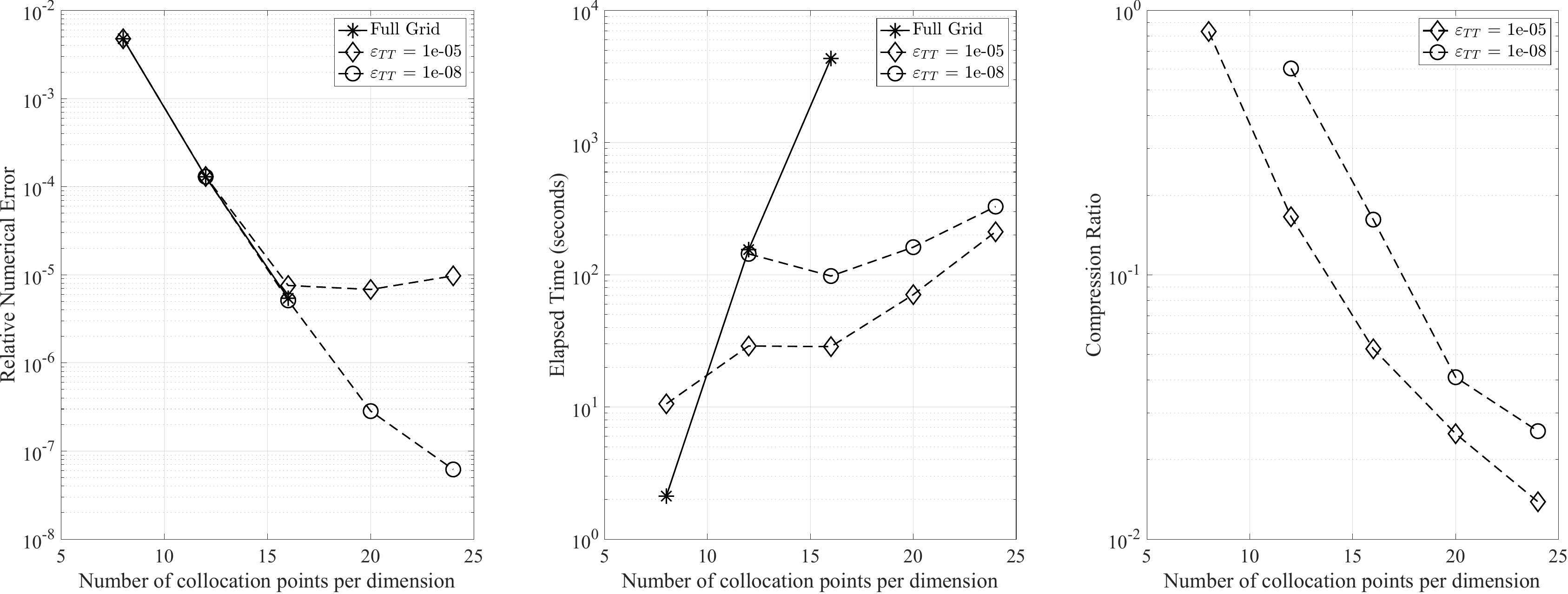}
    \caption{Experiment 3: Left panel - Relative numerical error versus number of collocation points per dimension. Middle panel - Elapsed time in seconds. Right panel - Compression ratio of the TT solutions.}
    \label{fig:burger}
\end{figure}

Figure~\ref{fig:burger} shows the performance of all three schemes. 
In the left panel, we plot the relative numerical error versus the number of collocation points per dimension. The TT schemes provide the same level of accuracy as the full grid scheme, and we see that the exponential convergence property of the spectral method is retained. Lastly, the error for TT schemes plateau as they reach their respective $\varepsilon_{TT}$.
In the middle panel, we show the elapsed computation time required for the three schemes. The result shows that the TT schemes are significantly more efficient than the full grid scheme. For example, at 16 collocation points per dimension, the TT scheme with $\varepsilon_{TT}=10^{-5}$ is about 150 faster than the full grid scheme. As the slope of the elapsed time curve for the full grid scheme is much steeper compared to TT schemes, the speedup will increase as the number of collocation points increases.
The right panel shows the compression ratio of the TT solutions for both TT schemes. The compression ratio becomes better as the number of collocation points increases. For 16 collocation points per dimension, the compression ratio is around $10^{-1}$, indicating the TT-ranks of the solutions are moderate.
For this realistic experiment, we find that TT scheme is more efficient than the full grid scheme, and potentially makes much higher resolution feasible.
\section{Conclusion}
In this work, we have developed a tensor network space-time spectral collocation method for solving the nonlinear convection-diffusion equation. Our approach leverages the Tensor Train (TT) format to mitigate the curse of dimensionality inherent in space-time discretizations. A key challenge is controlling the growth of TT-ranks during Newton iterations when solving the resulting nonlinear tensor equations. To address this, we introduced the "Step Truncation TT-Newton" algorithm that adapts the TT truncation tolerance at each iteration to constrain the solution to a low-rank tensor manifold.
Our numerical experiments demonstrated the exponential convergence of the proposed scheme through various benchmark problems. Importantly, the TT format provided a significantly reduced memory requirement compared to traditional full-grid discretizations. The step-truncation approach successfully maintained low-rank structures throughout Newton iterations, enabling substantial computational savings over full-grid solvers, especially at higher resolutions.

\section*{CRediT authorship contribution statement}
D. Adak, M. E. Danis, D. P. Truong, K. O. Rasmussen, B. S. Alexandrov : Conceptualization, Methodology, Writing-original draft, Review \& Editing.

\section*{Data availability statement}

The data that support the findings of this research are available from the corresponding author upon reasonable request.

\section*{ Declaration of competing interest}

The authors declare that they have no known competing financial interests or personal relationships that could have appeared to influence the work reported in this paper.


\section*{Acknowledgments}
The authors gratefully acknowledge the support of the Laboratory
Directed Research and Development (LDRD) program of Los Alamos
National Laboratory under project number 20230067DR.
Los Alamos National Laboratory is operated by Triad National Security,
LLC, for the National Nuclear Security Administration of
U.S. Department of Energy (Contract No.\ 89233218CNA000001).


\bibliographystyle{plain}
\bibliography{main}


\appendix
\section{}
\label{sec:app}
\subsection{Construction of Input Functions}
\label{app:func_eval}
At every Newton's iteration, input functions such as coefficient functions, loading tensor and the boundary term needs to be evaluated at the current solution $\tentt{U}$. 
Depending on their analytical forms, these functions can be evaluated either directly by available tensor train arithmetic \cite{oseledets2011tensor}, or by the cross interpolation. Subsequently, they need to be recompressed with a desired TT tolerance $\varepsilon_{TT}$ to achieve lower ranks if possible.
Next, after being evaluated, coefficient functions need to be converted into a TT-matrix operator for the multiplication with other operators. This procedure is the analogous TT version of converting a vector into a diagonal matrix. The algorithm to achieve this is described in \cite[Algorithm 1]{adak2024tensor}.

\subsection{Construction of the boundary term \texorpdfstring{$\ten{F}^{bc}(\tentt{U})$}{}} Here we provide details about the construction of $\ten{F}^{bc}(\tentt{U}) = \ten{S}^{map,TT}(\tentt{U})\ten{G}^{bc,TT}$.

\begin{equation*}
  \ten{S}^{map,TT} = \sum_{k = t,D,C} \ten{A}_k^{map,TT},
\end{equation*}
where 
\begin{equation}
\begin{split}
  &\ten{A}^{map,TT}_t = \bold{S}_t(\ten{I}_t,\ten{I}_t)\tenprod \bold{I}^{N+1}(\ten{I}_s,:)\tenprod \bold{I}^{N+1}(\ten{I}_s,:)\tenprod \bold{I}^{N+1}(\ten{I}_s,:),\\
\end{split}
\end{equation}
\begin{equation}
  \begin{split}
    &\ten{L}^{map,TT}=
    \begin{aligned}[t]
      &\bold {I}^{N+1}(\ten{I}_t,:)\tenprod \bold{S}_{xx}(\ten{I}_s,:) \tenprod \bold {I}^{N+1}(\ten{I}_s,:)\tenprod \bold {I}^{N+1}(\ten{I}_s,:)\\
      +&\bold {I}^{N+1}(\ten{I}_t,:) \tenprod \bold{I}^{N+1}(\ten{I}_s,:) \tenprod \bold{S}_{yy}(\ten{I}_s,:)  \tenprod\bold{I}^{N+1}(\ten{I}_s,:)\\
      +&\bold {I}^{N+1}(\ten{I}_t,:) \tenprod \bold{I}^{N+1(\ten{I}_s,:} \tenprod  \bold{I}^{N+1(\ten{I}_s,:} \tenprod\bold{S}_{zz}(\ten{I}_s,:),
    \end{aligned}\\
    &\ten{A}^{map,TT}_D(\tentt{U}) = \text{diag}(a(\tentt{U}))\ten{L}^{map,TT},
  \end{split}
\end{equation}

\begin{equation}
  \begin{split}
    &\nabla_x^{map,TT} = \bold {I}^{N+1}(\ten{I}_t,:)\tenprod \bold{S}_{x}(\ten{I}_s,:) \tenprod \bold{I}^{N+1}(\ten{I}_s,:)\tenprod \bold{I}^{N+1}(\ten{I}_s,:), \\
    &\nabla_y^{map,TT} = \bold {I}^{N+1}(\ten{I}_t,:)\tenprod \bold{I}^{N+1}(\ten{I}_s,:) \tenprod\bold{S}_{y}(\ten{I}_s,:) \tenprod \bold{I}^{N+1}(\ten{I}_s,:), \\
    &\nabla_z^{map,TT} = \bold {I}^{N+1}(\ten{I}_t,:)\tenprod \bold{I}^{N+1}(\ten{I}_s,:) \tenprod \bold{I}^{N+1}(\ten{I}_s,:) \tenprod\bold{S}_{z}(\ten{I}_s,:), \\
    &\ten{A}_C^{map,TT}(\tentt{U}) = \sum_{l=x,y,z} \text{diag}(b^l(\tentt{U}))\nabla_l^{map,TT}.
  \end{split}
\end{equation}

Next, we show how to construct the $\ten{G}^{bc,TT}$ tensor.
The tensor $\ten{G}^{Bd}$ is a $(N+1)\times (N+1) \times (N+1) \times
(N+1)$, in which only the BC/IC elements are computed. Other elements
are zeros. The TT tensor, $\ten{G}^{bc,TT}$ is constructed using the cross
interpolation.

Then, the boundary tensor $\ten{F}^{bc}(\tentt{U})$ is computed as:
\begin{equation}
  \ten{F}^{bc}(\tentt{U})= \ten{S}^{map,TT}(\tentt{U})\ten{G}^{bc,TT}
\end{equation}

\end{document}
